\documentclass[11pt]{article}%
\usepackage{amsmath}
\usepackage{amsfonts}
\usepackage{amssymb}
\usepackage{graphicx}
\usepackage{a4wide}%
\usepackage{hyperref}%
\newtheorem{theorem}{Theorem}
\newtheorem{example}[theorem]{Example}

\begin{document}

\title{On \textquotedblleft two important theorems\textquotedblright\ in canonical
duality theory}
\author{C. Z\u{a}linescu}
\date{}
\maketitle

\begin{abstract}
In this short note we show, providing counterexamples, that the
\textquotedblleft two important theorems\textquotedblright\ in the
recent paper [Y, Yuan, Global optimization solutions to a class of
non-convex quadratic minimization problems with quadratic
constraints, in Canonical Duality Theory, D.Y. Gao et al.\ (eds),
(AMMA, volume 37), Springer, 2017] are false.

\end{abstract}

\section{Introduction}

The aim of \cite{Yua:17} (the same as \cite{Yua:16}) is to study
\textquotedblleft non-convex quadratic minimization problems with quadratic
constraints\textquotedblright. The method used in that study is
\textquotedblleft a very powerful method, proposed by David Gao, called
\emph{canonical duality}\textquotedblright. It is our aim to show that both
theorems stated in this paper are false.

To ease the reading of our note we present the statements, as well
as their ingredients, making so the note self-contained. Of course,
we encourage the reader to look also at \cite{Yua:16} and/or
\cite{Yua:17}.

\section{Framework and statements of \cite{Yua:17}}

Let us quote the corresponding text from \cite{Yua:17}, split by some remarks.
So, the

\medskip\textquotedblleft Non-convex quadratic minimization problems with
quadratic constraints $(\mathcal{P}_{qq})$ can be formulated as follows
(($\mathcal{P}_{qq}$) in short)

$(\mathcal{P}_{qq})$ : $\min\left\{  P(x)=\tfrac{1}{2}x^{T}Ax-f^{T}%
x:x\in\mathcal{X}_{a}\right\}  ,\quad(1)$

\noindent where $A=A^{T}\in\mathbb{R}^{n\times n}$ is an indefinite matrix,
and the feasible space $\mathcal{X}_{a}$ is defined by

$\mathcal{X}_{a}\overset{\Delta}{=}\left\{  x\in\mathbb{R}^{n}\mid\tfrac{1}%
{2}x^{T}Q_{i}x+b_{i}^{T}x\leq c_{i},\ i=1,...,m\right\}  ,\quad(2)$

\noindent in which $Q_{i}=Q_{i}^{T}\in\mathbb{R}^{n\times n}$ $(i=1,...,m)$
are given nonsingular matrices, $b_{i}\in\mathbb{R}^{n}$ $(i=1,...,m)$ are
given vectors which control the geometric centers. $c_{i}(i=1,...,m)\in
\mathbb{R}$ are given input constants.

In order to make sure that the feasible space $\mathcal{X}_{a}$ is nonempty,
the quadratic constraints must satisfy the \emph{Slater regularity condition},
i.e., there exists one point $x_{0}$ such that $\tfrac{1}{2}x_{0}^{T}%
Q_{i}x_{0}+b_{i}^{T}x_{0}\leq c_{i}$, $i=1,...,m.$\textquotedblright

\medskip

First observe that the Slater regularity condition (in fact Slater's
constraint qualification, see \cite[p.\ 243]{BSS}) for problem $(\mathcal{P}%
_{qq})$ is: there exists $x_{0}\in\mathbb{R}^{n}$ such that $\tfrac{1}{2}%
x_{0}^{T}Q_{i}x_{0}+b_{i}^{T}x_{0}<c_{i}$ for all $i\in\overline{1,m}.$

\medskip

One continues with:

\medskip\textquotedblleft In this work, one hard restriction is given that
$f\neq0\in\mathbb{R}^{n}$. The restriction is very important to guarantee the
uniqueness of global optimal solution of $(\mathcal{P}_{qq})$%
\textquotedblright.

\medskip

In fact the condition $f\in\mathbb{R}^{n}\setminus\{0\}$ \textbf{does not
guarantee} the uniqueness of a global optimal solution of $(\mathcal{P}%
_{qq}),$ as Example \ref{ex1} (below) shows.

\medskip

The following notation is used in the statement of Theorem 2 of \cite{Yua:17}:

\medskip\textquotedblleft$\mathcal{G}_{+}(A)\overset{\Delta}{=}\{B\in
\mathbb{R}^{n\times n}\mid A+B\succ0\}.\quad(5)$\textquotedblright

\medskip

On page 344 of \cite{Yua:17} one finds:

\medskip\textquotedblleft The canonical dual function of $P(x)$ is defined by
the following equation (referred to [8])\footnote{The reference [8] is
\textquotedblleft Gao, DY. Canonical dual transformation method and
generalized triality theory in nonsmooth global optimization. J Global Optim
2000; 17(1--4): 127--160\textquotedblright.}

$P^{d}(\sigma)=Q^{\Lambda}(\sigma)-\mathcal{I}^{\ast}(\sigma),\quad(10)$

\noindent where

$Q^{\Lambda}(\sigma)=sta\left\{  \varepsilon^{T}\sigma+\tfrac{1}{2}%
x^{T}Ax-f^{T}x\right\}  =-\tfrac{1}{2}F(\sigma)^{T}G(\sigma)^{-1}%
F(\sigma)-c^{T}\sigma,\quad(11)$

\noindent in which the notation $sta\left\{  ^{\ast}:x\in\mathbb{R}%
^{n}\right\}  $ is the operator to find out the stationary point in the space
$\mathbb{R}^{n}$, $G(\sigma),$ $F(\sigma)$ and $c$ are defined by

$G(\sigma)=\left(  A+\sum_{i=1}^{m}Q_{i}\sigma_{i}\right)  ,$ $F(\sigma
)=\left(  f-\sum_{i=1}^{m}b_{i}\sigma_{i}\right)  ,$ $c=(c_{1},c_{2}%
,\ldots,c_{m})^{T},\quad(12)$

\noindent where $\sigma_{i}$ is the $i$th element of $\sigma$.

The dual feasible space is defined by

$\mathcal{S}\overset{\Delta}{=}\left\{  \sigma\in\mathbb{R}^{m}\mid\sigma
\geq0\in\mathbb{R}^{m},\ \det(G(\sigma))\neq0\right\}  .\quad(13)$%
\textquotedblright.

\medskip

Above, $\mathcal{I}^{\ast}$ is given by $\mathcal{I}^{\ast}(\sigma)=0$ if
$\sigma\geq0,$ $\mathcal{I}^{\ast}(\sigma)=+\infty$ otherwise, and so

\medskip$P^{d}(\sigma)=-\tfrac{1}{2}F(\sigma)^{T}G(\sigma)^{-1}F(\sigma
)-c^{T}\sigma$ if $\sigma\geq0,$ $P^{d}(\sigma)=-\infty$ otherwise.

\medskip Notice that $P^{d}(\sigma)$ is not defined if $\sigma\geq0$ and
$\det(G(\sigma))=0!$ Also, it is quite strange that an \textquotedblleft
operator that is used to find the stationary point in the space $\mathbb{R}%
^{n}$\textquotedblright\ could be a real number.

\medskip

One continues with:

\medskip\textquotedblleft The canonical dual problem ($\mathcal{P}^{d}$ in
short) associated with $(\mathcal{P}_{qq})$ can be eventually formulated as follows

$(\mathcal{P}^{d})$ : $\max_{\sigma\in\mathcal{S}}\left\{  P^{d}%
(\sigma)\right\}  .\quad(14)$

\medskip

\textbf{2.3 Two important theorems}

In order to show that there is no duality gap, the following theorem is presented.

\textbf{Theorem 1.} If $A,Q_{i},b_{i},f_{i},c_{i},$ $i=1,2,.,m$, are given
with definitions in $(\mathcal{P}_{qq})$ such that the dual feasible space

$\mathcal{Y}\overset{\Delta}{=}\left\{  \sigma\in\mathcal{S}\mid
G(\sigma)^{-1}F(\sigma)\in\mathcal{X}\right\}  \quad(15)$

\noindent is not empty, then the problem

$(\mathcal{P}^{d})$ : $\max_{\sigma\in\mathcal{Y}}\left\{  P^{d}%
(\sigma)\right\}  ,\quad(16)$

\noindent is canonically (perfectly) dual to $(\mathcal{P}_{qq})$. In another
words, if $\overline{\sigma}$ is a solution of the dual problem $(\mathcal{P}%
^{d})$,

$\overline{x}=G(\overline{\sigma})^{-1}F(\overline{\sigma})\quad(17)$

\noindent is a solution of $(\mathcal{P}_{qq})$ and

$P(\overline{x})=P^{d}(\overline{\sigma}).\quad(18)$\textquotedblright

\medskip

After the proof of this \textquotedblleft important theorem\textquotedblright,
one continues with:

\medskip\textquotedblleft In order to get the optimization solution of
$(\mathcal{P}_{qq} )$, we introduce the following subset

$\mathcal{S}_{+}=\left\{  \sigma\in\mathcal{S}\mid G(\sigma)\text{ is positive
definite}\right\}  .\quad(23)$

In order to hold on the uniqueness of the optimal duality solution, the
following existence theorem is presented.

\textbf{Theorem 2.} For any given symmetrical matrixes $A,Q_{i},\in
\mathbb{R}^{n\times n}$, $\mathcal{G}_{+}(A)$ (defined by (5)) is the
complementary positive definite matrix group of $A$, $f$, $b_{i}\in
\mathbb{R}^{n},$ $c_{i}\in\mathbb{R},$ $i=1,2,...,m,$ if the following two
conditions are satisfied

$C_{1}:\sum_{i=1}^{m}Q_{i}\in\mathcal{G}_{+}(A);$

$C_{2}:$ there must exist one $k(1\leq k\leq m)$ such that $Q_{k}$ is
positive-definite and $Q_{k}\in\mathcal{G}_{+}(A)$, moreover,

$\left\Vert D_{k}A^{-1}f\right\Vert >\left\Vert b_{k}^{T}D_{k}^{-1}\right\Vert
+\sqrt{\left\Vert b_{k}^{T}D_{k}^{-1}\right\Vert ^{2}+2\left\vert
c_{k}\right\vert },\quad(24)$

\noindent where $Q_{k}=D_{k}^{T}D_{k}$ and $\left\Vert ^{\ast}\right\Vert $ is
some vector norm.

Then, the canonical duality problem (16) has a unique nonzero solution
$\overline{\sigma}$ in the space $\mathcal{S}_{+}$.\textquotedblright%
\footnote{Observe that $A$ is not assumed to be non-singular!}

\medskip

First of all observe that there are two dual problems, (14) and (16). However,
both theorems refer to problem (16). Moreover, probably the author intended to
write $\mathcal{X} _{a}$ instead of $\mathcal{X}$, and $f$ instead of $f_{i}$
in the statement of Theorem 1. Because $Q_{k}$ is positive-definite in Theorem
2, $D_{k}:=Q_{k}^{1/2}$ is positive-definite (hence symmetric) and $D_{k}%
D_{k}=Q_{k};$ moreover, the norm has to be the Euclidean norm to have correct
inequalities in the proof of Theorem 2. Furthermore, I suppose that the notion
of solution of the problem $(\mathcal{P}_{qq})$ is in the sense from
\cite[p.\ 2]{BSS}; similarly for solution of problem $(\mathcal{P}^{d}).$

\section{Examples}

The first example shows that the condition $f\neq0\in\mathbb{R}^{n}$ does not
\textquotedblleft guarantee the uniqueness of global optimal solution of
$(\mathcal{P}_{qq})$\textquotedblright.

\begin{example}
\label{ex1}Take $n=m=1,$ $P(x):=-x^{2}+2x$ and $q(x):=-P(x)$ for
$x\in\mathbb{R}.$ Clearly, the problem $\min$ $P(x)$ s.t. $q(x)\leq0$ has the
solutions $x_{1}=0$ and $x_{2}=2.$ In fact taking an arbitrary quadratic
function $P$ on $\mathbb{R}^{n}$ and $q:=-P,$ assuming that $\mathcal{X}%
_{a}:=\{x\in\mathbb{R}^{n}\mid q(x)\leq0\}\neq\emptyset,$ the sets of
solutions of the problem $\min$ $P(x)$ s.t. $q(x)\leq0$ is $\{x\in
\mathbb{R}^{n}\mid q(x)=0\}$.
\end{example}

The next example shows that $\overline{x}$ provided in \cite[Eq.\ (17)]%
{Yua:17} is not, necessarily, a solution of the primal problem \cite[Eq.\ (1)]%
{Yua:17}, contrary to what is stated in \cite[Th.\ 1]{Yua:17}.

\begin{example}
\label{ex2}Let us take $n=m=1,$ $P(x):=-x^{2}-x$ and $q(x):=\tfrac{1}{2}%
x^{2}+x$ and $c:=0;$ hence $\mathcal{X}_{a}=[-2,0].$ With the notations in
\cite[Eq.~(12)]{Yua:17} we have that $G(\sigma)=-2+\sigma,$ $F(\sigma
)=1-\sigma,$ and so $P^{d}(\sigma)=-\tfrac{1}{2}\frac{(1-\sigma)^{2}}%
{\sigma-2}$ (for $\sigma\in\mathcal{S}=[0,2)\cup(2,\infty)$). We have that
$\mathcal{Y}=\{\sigma\in\mathcal{S}\mid q(\frac{1-\sigma}{\sigma-2}%
)\leq0\}=[0,1]\cup\lbrack3,\infty).$ It is easy to verify that $P^{d}(0)\geq
P^{d}(\sigma)$ for all $\sigma\in\mathcal{Y}.$ In fact $\overline{\sigma}:=0$
is the unique solution of problem $(\mathcal{P}^{d})$. Indeed, we have that
$\overline{x}=F(\overline{\sigma})/G(\overline{\sigma})=-\tfrac{1}{2}.$
However, $\tfrac{1}{4}=P^{d}(\overline{\sigma})=P(\overline{x})>P(-2)=-2.$
Hence $\overline{x}$ is not a solution of problem $(\mathcal{P}_{qq});$ in
fact $P(\overline{x})>P(x)$ for every $x\in\mathbb{R}\setminus\{\overline
{x}\}.$
\end{example}

The next example shows that under the hypothesis of \cite[Th.\ 2]{Yua:17}, its
conclusion that the dual problem \cite[Eq.\ (16)]{Yua:17} has a unique
solution in $\mathcal{S}_{+}$ can be false; in fact in this example the dual
problem \cite[Eq.\ (16)]{Yua:17} has no solutions belonging to the set
$\mathcal{S}_{+}.$

\begin{example}
Let $A:=\left(
\begin{array}
[c]{ll}%
1 & 0\\
0 & -1
\end{array}
\right)  ,$ $Q:=\left(
\begin{array}
[c]{ll}%
4 & 0\\
0 & 4
\end{array}
\right)  \succ0$, $f:=\left(
\begin{array}
[c]{c}%
\sqrt{27}\\
1
\end{array}
\right)  ,$ $b:=\left(
\begin{array}
[c]{c}%
0\\
0
\end{array}
\right)  ,$ and $c=52.$ Clearly, $A+Q\succ0,$ and so condition $C_{1}$ holds.
Moreover, $D^{2}=Q,$ where $D=\tfrac{1}{2}Q,$ and $\left\Vert DA^{-1}%
f\right\Vert =\left\Vert (2\sqrt{27},-2)\right\Vert =4\sqrt{7}>\sqrt
{2c}=2\sqrt{26}.$ Hence, \cite[Eq.~(24)]{Yua:17} is verified, and so condition
$C_{2}$ holds, too. The arguments for the claim that \cite[Eq.\ (16)]{Yua:17}
has no solutions belonging to the set $\mathcal{S}_{+}$ are provided below.
\end{example}

The inequality constraint is $q(x_{1},x_{2})=\tfrac{1}{2}\left(  4x_{1}%
^{2}+4x_{2}^{2}\right)  -52\leq0$

With the notations in \cite[Eq.~(12)]{Yua:17} we have that $G(\sigma
)=\operatorname*{diag}(4\sigma+1,4\sigma-1),$ $F(\sigma)=f,$ $c=52,$ and so
$P^{d}(\sigma)=-\tfrac{1}{2}\left(  \frac{27}{4\sigma+1}+\frac{1}{4\sigma
-1}\right)  -52\sigma$ (for $\sigma\in\mathcal{S}=[0,\tfrac{1}{4})\cup
(\tfrac{1}{4},\infty)$). Moreover, $\mathcal{S}_{+}=(\tfrac{1}{4},\infty).$
Observe that
\[
\psi(\sigma):=q(G(\sigma)^{-1}F(\sigma))=q\left(  \frac{\sqrt{27}}{4\sigma
+1},\frac{1}{4\sigma-1}\right)  =2\left[  \frac{27}{\left(  4\sigma+1\right)
^{2}}+\frac{1}{\left(  4\sigma-1\right)  ^{2}}-26\right]
\]
for $\sigma\in\mathbb{R}\setminus\left\{  -\tfrac{1}{4},\tfrac{1}{4}\right\}
.$ Let us study the functions $\psi$ and $\varphi$ defined by
\[
\varphi:\mathbb{R}\setminus\left\{  -\tfrac{1}{4},\tfrac{1}{4}\right\}
\rightarrow\mathbb{R},\quad\varphi(\sigma):=-\tfrac{1}{2}\left(  \frac
{27}{4\sigma+1}+\frac{1}{4\sigma-1}\right)  -52\sigma.
\]
Clearly
\[
\varphi^{\prime}(\sigma)=2\left[  \frac{27}{\left(  4\sigma+1\right)  ^{2}%
}+\frac{1}{\left(  4\sigma-1\right)  ^{2}}-26\right]  =\psi(\sigma
),\quad\varphi^{\prime\prime}(\sigma)=-16\left[  \frac{27}{\left(
4\sigma+1\right)  ^{3}}+\frac{1}{\left(  4\sigma-1\right)  ^{3}}\right]
\]
for $\sigma\in\mathbb{R}\setminus\left\{  -\tfrac{1}{4},\tfrac{1}{4}\right\}
.$ First, we have the following table of variation for $\psi=\varphi^{\prime
}:$

\medskip

\begin{center}
{\footnotesize
\begin{tabular}
[c]{c|ccccccccccccccccccc}%
$\sigma$ & $-\infty$ &  & $\sigma_{1}$ &  & $-\tfrac{1}{4}$ &  & $0$ &  &
$\sigma_{2}$ &  & $\tfrac{1}{8}$ &  & $\sigma_{3}$ &  & $\tfrac{1}{4}$ &  &
$\sigma_{4}$ &  & $\infty$\\\hline
$\psi^{\prime}$ & $0$ & $+$ & $+$ & $+$ & $|$ & $-$ & $-$ & $-$ & $-$ & $-$ &
$0$ & $+$ & $+$ & $+$ & $|$ & $-$ & $-$ & $-$ & $0$\\\hline
$\psi$ & $-$ & $\nearrow$ & $0$ & $\nearrow$ & $|$ & $\searrow$ & $+$ &
$\searrow$ & $0$ & $\searrow$ & $-$ & $\nearrow$ & $0$ & $\nearrow$ & $|$ &
$\searrow$ & $0$ & $\searrow$ & $-$\\\hline
$\psi^{\prime\prime}$ &  & $+$ & $+$ & $+$ & $|$ & $+$ & $+$ & $+$ & $+$ & $+$
& $+$ & $+$ & $+$ & $+$ & $|$ & $+$ & $+$ & $+$ &
\end{tabular}
},
\end{center}

\medskip\noindent where we have taken into consideration that $\psi^{\prime
}(1/8)=0,$ $\psi(1/8)=-20<0,$ and%
\[
\lim_{\sigma\rightarrow\pm1/4}\left\vert \psi^{\prime}(\sigma)\right\vert
=\lim_{\sigma\rightarrow\pm1/4}\psi(\sigma)=\infty,\quad\lim_{\sigma
\rightarrow\pm\infty}\psi^{\prime}(\sigma)=0,\quad\lim_{\sigma\rightarrow
\pm\infty}\psi(\sigma)=-52.
\]

Taking into account the variation of $\varphi^{\prime}=\psi,$ we get the
following table for the variation of $\varphi$ on the interval $[0,\infty):$

\medskip

\begin{center}
{%
\begin{tabular}
[c]{c|ccccccccccccc}%
$\sigma$ & $0$ &  & $\sigma_{2}$ &  & $\tfrac{1}{8}$ &  & $\sigma_{3}$ &  &
$\tfrac{1}{4}$ &  & $\sigma_{4}$ &  & $\infty$\\\hline
$\varphi^{\prime}$ & $4$ & $+$ & $0$ & $-$ & $-$ & $-$ & $0$ & $+$ &
$^{+\infty}|^{+\infty}$ & $+$ & $0$ & $-$ & $-$\\\hline
$\varphi$ & $-13$ & $\nearrow$ &  & $\searrow$ & $-$ & $\searrow$ &  &
$\nearrow$ & $^{+\infty}|_{-\infty}$ & $\nearrow$ &  & $\searrow$ &
$-$\\\hline
$\varphi^{\prime\prime}$ & $-$ & $-$ & $-$ & $-$ & $0$ & $+$ & $+$ & $+$ & $|$
& $-$ & $-$ & $-$ &
\end{tabular}
.}
\end{center}

\medskip

We obtain that $\mathcal{X}_{a}=\{x\in\mathbb{R}^{2}\mid\left\Vert
x\right\Vert \leq\sqrt{26}\},$ $\mathcal{S}=[0,\tfrac{1}{4})\cup(\tfrac{1}%
{4},\infty),$ $\mathcal{Y}=[\sigma_{2},\sigma_{3}]\cup\lbrack\sigma_{4}%
,\infty).$ Since $\varphi$ is decreasing on each one of the intervals
$[\sigma_{2},\sigma_{3}]$ and $[\sigma_{4},\infty),$ $\max_{\sigma
\in\mathcal{Y}}P^{d}(\sigma)=\max\{\varphi(\sigma_{2}),\varphi(\sigma_{4})\}.$
Clearly, $\varphi(\sigma_{2})\geq\varphi(\tfrac{1}{8})=-29/2.$ So, it is
sufficient to prove that $\varphi(\sigma)\leq-29/2$ on $(\tfrac{1}{4}%
,\infty).$ Replacing $4\sigma$ by $t,$ the preceding inequality is equivalent
to each one of the following inequalities $\frac{27}{t+1}+\frac{1}%
{t-1}+26t\geq29,$ $\chi(t):=26t^{3}-29t^{2}+2t+3\geq0$ for $t>1.$ But
$\chi^{\prime}(t)=78t^{2}-58t+2=27t(39t-29)+2>0$ for $t\geq1,$ whence
$\chi(t)\geq\chi(1)=2$ for $t\geq1.$ Hence $\sigma_{2}$ $(\in\mathcal{Y}%
\setminus\mathcal{S}_{+})$ is the only solution of problem (16). This shows
that the conclusion of Theorem 2 in \cite{Yua:17} is false.

\end{document}